\documentclass[a4paper,12pt]{amsart}

\usepackage[centertags]{amsmath}
\usepackage{amsfonts}
\usepackage{amssymb}
\usepackage{amsthm}
\usepackage{newlfont}
\usepackage{euscript}

\newtheorem{thm}{Theorem}
\newtheorem{lem}{Lemma}
\newtheorem{prop}{Proposition}
\theoremstyle{remark}
\newtheorem{zau}{Remark}
\theoremstyle{definition}
\newtheorem{ozn}{Definition}
\newtheorem{ex}{Example}
\begin{document}

\newcommand{\demo}{\emph{Proof.} }

\newcommand{\1}{1\!\!{\mathrm I}}
\renewcommand{\Re}{{\Bbb R}}
\newcommand{\eps}{\varepsilon}
\newcommand{\kap}{\varkappa}
\newcommand{\vO}{\varOmega}
\newcommand{\bu}{\mathbf{u}}
\newcommand{\bp}{\mathbf{p}}
\newcommand{\ax}{\Re^+}
\newcommand{\prt}{\partial}
\newcommand{\se}{\mathsf{E}}
\newcommand{\cs}{\mathsf{C}}
\newcommand{\Ff}{{\EuScript F}}
\newcommand{\Sf}{\mathcal{S}}
\newcommand{\Ef}{{\EuScript E}}
\newcommand{\Bf}{{\EuScript B}}
\newcommand{\Tf}{{\EuScript T}}
\newcommand{\Cf}{{\EuScript C}}
\newcommand{\Zf}{{\EuScript Z}}
\newcommand{\Wf}{{\EuScript W}}
\newcommand{\Df}{{\EuScript D}}
\newcommand{\Nf}{{\EuScript N}}
\newcommand{\Gf}{{\EuScript G}}
\newcommand{\Qf}{{\EuScript Q}}
\newcommand{\Lf}{{\EuScript L}}
\newcommand{\Uf}{{\EuScript U}}
\newcommand{\pf}{\Pi_{fin}}
\newcommand{\Kb}{\mathbf{K}}
\newcommand{\ZZ}{{\Bbb Z}}
\newcommand{\NN}{{\Bbb N}}
\newcommand{\QQ}{{\Bbb Q}}
\newcommand{\TT}{{\Bbb T}}
\newcommand{\GG}{{\Bbb G}}
\newcommand{\DD}{{\Bbb D}}
\newcommand{\NNN}{\NN\times \NN}
\newcommand{\Cd}{C_\bullet\,}
\newcommand{\Span}{\mathrm{span}\,}
\newcommand{\bro}{\hbox{{\boldmath $\rho$}}}
\newcommand{\sbro}{\hbox{{\boldmath $\rho$}}}
\newcommand{\brosup}{\hbox{{\boldmath $\rho_{\sup}$}}}
\newcommand{\bkap}{\hbox{{\boldmath $\vartheta$}}}
\newcommand{\bnu}{\hbox{{\boldmath $\nu$}}}
\newcommand{\bpi}{\hbox{{\boldmath $\Pi$}}}
\newcommand{\sbpi}{\hbox{{\small\boldmath $\Pi$}}}
\newcommand{\bg}{\hbox{{\boldmath $\Gamma$}}}
\newcommand{\ba}{\hbox{{\boldmath $\alpha$}}}
\newcommand{\sba}{\hbox{{\boldmath $\alpha$}}}
\newcommand{\sbg}{\hbox{{\small\boldmath $\Gamma$}}}
\newcommand{\pa}{\prt_\alpha}
\newcommand{\pba}{\prt_{\sba}}
\newcommand{\eqdef}{\mathop{=}\limits^{df}}
\newcommand{\be}{\begin{equation}}
\newcommand{\ee}{\end{equation}}


\title[Distribution density for an Ornstein-Uhlenbeck
process with L\'evy noise] {Conditions for existence and
smoothness of the distribution density for an Ornstein-Uhlenbeck
process with L\'evy noise}
\author{Semen V.Bodnarchuk, Alexey M.Kulik}
\address{Kiev 03022 Glushkova str. 6, Kyiv Taras Shevchenko National University, Ukraine} \email{sem\_bodn@ukr.net}

\address{Kiev 01601 Tereshchenkivska str. 3, Institute of Mathematics,
Ukrai\-ni\-an National Academy of Sciences, Ukraine}
 \email{kulik@imath.kiev.ua}

\begin{abstract} Conditions are given, sufficient for the distribution of
an Ornstein-Uhlenbeck process with L\'evy noise to be absolutely
continuous or to possess a smooth density. For the processes with
non-degenerate drift coefficient, these conditions are a necessary
ones. A multidimensional analogue for the non-degeneracy condition
on the drift coefficient is introduced.
\end{abstract}


\maketitle


\section{Introduction}
The theory of stochastic differential equations (SDE) with jump
noise is developed intensively during the last decades. It is
stimulated by a wide range of the disciplines that use such types
of SDE's as a models,  from climatology (e.g.  \cite{im}) to
financial mathematics (e.g. \cite{cont}, \cite{geman}). One of the
most important point in this theory consists in studying a local
properties of the laws of the solutions to such an equations. For
instance, an information about the distribution density for the
solution allows one to investigate effectively the ergodic
properties of this solution  (see \cite{Me_ergodic} and discussion
therein). This, in turn, allows one  to conduct consistent
statistical analysis for such a processes, to solve a filtration
and optimal control problems for such a processes, etc.

The large variety of publications is devoted to investigation of
the properties of the laws of solutions to SDE's with jump noise
(e.g. \cite{bismut} -- \cite{Nou_sim}). These properties depend
essentially both on the structure of the equation and its
coefficients, and on the characteristics of the jump noise (i.e.,
its L\'evy measure). Although a wide spectrum of a sufficient
conditions is available, these conditions are not completely
satisfactory and can not be considered as a definitive ones. On
the one hand, it is difficult to compare the available sufficient
conditions. On the other hand, it is unclear how close   these
conditions are to the necessary ones. Therefore, an important (and
non simple) question is about the proper form of sufficient
conditions for existence and smoothness of the distribution
density, close to the necessary ones. In this article, we give an
answer to this question in the structurally most simple class of
linear SDE's with additive jump noise. Solutions to such an
equations often are called Ornstein-Uhlenbeck processes with
L\'evy noise.

\section{Formulation of the problem}

Consider linear SDE  in $\Re^m$,
\begin{equation}\label{001}
X(t)=X(0)+\int_0^tAX(s)\,ds+Z(t),
\end{equation}
where $X(0)\in\Re^m, A$ is an $m\times m$-matrix, $Z$ is an
$\Re^m$-valued L\'evy process (i.e., a continuous in probability
time homogeneous process with independent increments). It is well
known (e.g. \cite{Skor_nez_prir})  that every such a process
possesses representation
\be\label{0011}Z(t)=Z(0)+at+BW(t)+\int_0^t\int_{\|u\|_{\Re^m}>1}u\nu(ds,du)
+\int_0^t\int_{\|u\|_{\Re^m}\leq 1}u\tilde\nu(ds,du), \ee where
$a\in\Re^m,B\in\Re^{m\times m}$ are  deterministic vector and
matrix, respectively, $W$ is the Wiener process in $\Re^m$, $\nu$
is the Poisson random point measure on $\ax\times \Re^m$ with its
intensity measure equal  $dt\times \Pi(du)$ ($\Pi$ is the L\'evy
measure of the measure $\nu$), and $\tilde
\nu(ds,du)=\nu(ds,du)-ds\Pi(du)$ is the corresponding compensated
measure, it being known that $W$ and $\nu$ are independent.

Equation (\ref{001}) can be naturally interpreted as a family of a
Volterra type integral equations, indexed by the probability
variable $\omega$. Thus, for every measurable process $Z$ with its
trajectories being a.e. locally bounded, this equation possesses
unique solution with its trajectories also being a.e. locally
bounded. Remark that every L\'evy process possesses a modification
that satisfies conditions on the process $Z$ formulated before,
and therefore the solution to (\ref{001}) is well defined.
Moreover, this solution has the explicit representation
\be\label{13}X(t)=e^{\,tA}X(0)+\int_0^te^{\,(t-s)A}a\,ds+\int_0^te^{\,(t-s)A}B\,dW(s)+
$$
$$
+\int_0^t\int_{\|u\|_{\Re^m}>1}e^{\,(t-s)A}u\nu(ds,du)
+\int_0^t\int_{\|u\|_{\Re^m}\leq
1}e^{\,(t-s)A}u\tilde\nu(ds,du),\quad  t\geq0,\ee where
$e^{tA}=\sum_{n=0}^\infty{(tA)^k\over k!}, t\geq 1$ is the
solution to the matrix-valued differential equation $dE(t)=AE(t)\,
dt, E(0)=I_{\Re^m}$ is the unit matrix $m\times m$. Formula
(\ref{13}) is verified straightforwardly via the Ito formula.

All the summands in (\ref{13}) are independent and the fourth one
takes value 0 on the set
$\{\nu([0,t]\times\{\|u\|_{\Re^m}>1\})=0\}$, the latter having
positive probability. Thus $X(t)$ possesses (smooth) distribution
density iff so does the sum in which this summand is absent.
Moreover, the first two summands in (\ref{13}) are deterministic
and obviously do not have effect on existence or smoothness of the
density. Thus we suppose in a sequel that $X(0)=0, a=0,
\Pi(\|u\|>1)=0$.

\section{One-dimensional equation}

In the one-dimensional case, $A,B$ are a real numbers and the
third summand in (\ref{13}) is a normal random variable with its
second moment equal $B^2\int_0^te^{2(t-s)A}\,ds$. Since the
summands in (\ref{13}) are independent, for $B\not=0$ the law of
$X(t)$ is the convolution of some distribution with a
non-degenerate Gaussian one, and thus possesses a smooth density.
Further in this section, we consider the case $B=0$.

In a separate case $A=0$ we have  $X(t)=Z(t)-Z(0)$, and the
question on the properties of the law of the solution to
(\ref{001}) is exactly the same question for the distribution of
the L\'evy process that does not contain a diffusion component.
The complete answer to this question is not available now. Let us
formulate two sufficient conditions.

\begin{prop} 1. \cite{sato82} Denote $\mu(du)=[u^2\wedge 1] \Pi(du)$.
If $\Pi(\Re)=+\infty$ and, for some $n\in\NN$, the $n$-th
convolution power of the measure $\mu$ is absolutely continuous,
then  the distribution of $Z(t)$  is absolutely continuous for
every $t>0$.

2. \cite{kallenberg} If
$\left[\varepsilon^2\ln\frac{1}{\varepsilon}\right]^{-1}\int_{-\eps}^\eps
u^2\Pi(du)\to +\infty, \eps\to 0+,$ then the distribution of
$Z(t)$ possesses the  $C_b^\infty$ density for every $t>0$.
\end{prop}

Here and below, we denote by $C_b^\infty$ the class of infinitely
differentiable functions, bounded with all their derivatives. In a
sequel, we refer to the conditions formulated in the parts  1 and
 2 of Proposition 1 as for the Sato condition and the Kallenberg condition, respectively.
We emphasize once more, that both these conditions are sufficient
ones, but none of them is necessary. It appears that, for the
equation with its drift coefficient being \emph{non-degenerate},
the necessary and sufficient conditions are available both for
existence of the distribution density and for smoothness of this
density.

\begin{prop} Let $B=0, A\not=0$. Then the law of $X(t)$ is absolutely continuous
for every $t>0$ iff $\Pi(\Re)=+\infty$.
\end{prop}

It is obvious, that the condition $\Pi(\Re)=+\infty$ is necessary:
if $\Pi(\Re)=Q<+\infty,$ then the law of $X(t)$ has an atom with
its mass equal $e^{-tQ}$. Sufficiency follows from more general
Theorem 4.3 \cite{Me_TViMc} or Theorem A \cite{Nou_sim}.

\begin{thm} The following three statements are equivalent:

(i) for every $t>0$, the variable $X(t)$ possesses a distribution
density from the class  $C_b^\infty$;

(ii) for every $t>0$, the variable $X(t)$ possesses a bounded
distribution density;

(iii)
$\left[\varepsilon^2\ln\frac{1}{\varepsilon}\right]^{-1}\int_{\Re}
(u^2\wedge \eps^2)\Pi(du)\to +\infty, \eps\to 0+.$
\end{thm}
\demo The implication (i) $\Rightarrow$ (ii) is obvious. Let us
prove  that (ii) $\Rightarrow$ (iii). Denote  $\rho(\eps)=
\left[\varepsilon^2\ln\frac{1}{\varepsilon}\right]^{-1}\int_{\Re}
(u^2\wedge \eps^2)\Pi(du)$. Take $\eps\in (0,1)$ and write
$$X(t)=\int_0^t\int_{|u|\leq \eps}e^{(t-s)A}u\tilde
\nu(ds,du)+\int_0^t\int_{|u|\in(\eps,1]}e^{(t-s)A}u\tilde
\nu(ds,du).$$ The second moment  of the first summand is estimated
by $te^{2|A|t}\int_{-\eps}^\eps u^2\Pi(du).$ Thus the Chebyshev
inequality yields that the probability for the modulus of this
summand not to exceed $\sqrt{\eps}$ is not less than
$$
1-\eps^{-1}te^{2|A|t}\int_{-\eps}^\eps u^2\Pi(du)\geq
1-te^{2|A|t}\left[\eps\ln{1\over \eps}\right]\rho(\eps).
$$
The second summand is equal
$M(t,\eps)=\int_0^t\int_{|u|\in(\eps,1]}e^{(t-s)A}u\Pi(du)ds$ with
probability not less than $P(\nu((0,t)\times
\{|u|\in(\eps,1]\})=0)=\exp[-t\Pi(|u|\in(\eps,1])]$. The latter
term is not less than  $\exp[-t\ln{1\over
\eps}\rho(\eps)]=\eps^{t\rho(\eps)}.$ Therefore, \be\label{21}
P(X(t)\in [M(t,\eps)-\sqrt\eps, M(t,\eps)+\sqrt{\eps}])\geq
\eps^{t\rho(\eps)}-te^{2|A|t}\left[\eps\ln{1\over
\eps}\right]\rho(\eps). \ee Let (iii) fail, then there exists a
sequence $\eps_n\to 0+$ such that $\rho(\eps_n)\leq C<+\infty$.
Then, for  $t<{1\over
2C},x_n=M(t,\eps_n)-\sqrt\eps_n,y_n=M(t,\eps_n)+\sqrt\eps_n$, it
follows from (\ref{21}) that \be\label{24}
{P(X(t)\in[x_n,y_n])\over y_n-x_n}\to +\infty, \quad n\to +\infty.
\ee This, in turn, implies that (ii)  fails, also.

Let us prove  that (iii) $\Rightarrow$ (i). We remark that, by
general properties of the Fourier transform, the following
condition on the characteristic function $\phi$ of a random vector
in $\Re^m$ is sufficient for this vector to possess a $C_b^\infty$
distribution density: \be\label{22} \forall \, n \geq0\quad
\|z\|^n_{\Re^m} |\phi(z)|\to 0, \quad \|z\|_{\Re^m}\to~\infty.\ee
This condition is well known; sometimes, it is called a
(C)-condition (e.g. \cite{kallenberg}).

The value $X(t)$ is an integral of a deterministic function over
the compensated Poisson point measure. Therefore, its
characteristic function can be expressed explicitly:
\begin{equation}\label{23}
\phi_{X(t)}(z)=\exp\left\{\int_0^t\int_\Re
\left[\exp\{ize^{(t-s)A}u\}-1-ize^{(t-s)A}u\right]
\Pi(du)ds\right\}.
\end{equation}
 Without loss of generality, we can suppose that $A>0$. In what follows, we suppose $t>0$
 to be fixed. Take  $\beta>0$ in such a way that $\beta e^{(t-s)A}\leq1, s\in[0,t]$
 (i.e.,  $\beta=e^{-At}$). Denote
$$I_1(s,z)=\int_{\{|uz|\leq
\beta\}}\left[\cos{(e^{(t-s)A}uz)}-1\right]\Pi(du),$$
$$
I_2(s,z)=\int_{\{|uz|>
\beta\}}\left[\cos{(e^{(t-s)A}uz)}-1\right]\Pi(du).$$ Then, by
(\ref{23}),
\begin{eqnarray*}
|\phi_{X(t)}(z)|=\exp\left\{\int_0^tI_1(s,z)\,
ds+\int_0^tI_2(s,z)\right\},\quad z\in\Re.
\end{eqnarray*}
Let $C={1-\cos 1}$, one can verify  that  $\cos x-1\leq-Cx^2,
|x|\leq 1$. Then
$$I_1(s,z)\leq-C\int_{|uz|\leq\beta}(e^{(t-s)A}uz)^2\, \Pi(du)=
-Cz^2e^{2(t-s)A}\int_{|uz|\leq\beta}u^2\, \Pi(du),$$ and thus
$\int_0^tI_1(s,z)\,ds\leq-C_1z^2\int_{|uz|\leq\beta}u^2\, \Pi(du)$
with $C_1=C\,\frac{e^{\,2tA}-1}{2A}>0$. Next,
$$\int_0^tI_2(s,z)\,ds=\int_{|uz|>\beta}\int_0^t
\left[\cos{(e^{(t-s)A}uz)}-1\right]\,ds\Pi(du).
$$
By making the change of the variables $s\mapsto y=e^{(t-s)A}uz$
and taking into account that the function $x\mapsto \cos x-1$ is
an even one, we get
$$\int_0^tI_2(s,z)\,ds=\frac{1}{A}\int_{|uz|>\beta}
\int_{uz}^{e^{\,tA}uz}\frac{\cos y-1}{y}\,dy\Pi(du)\leq
$$
$$
\leq\int_{|uz|>\beta}\frac{1}{A|uz|}
\int_{|uz|}^{e^{\,tA}|uz|}(\cos y-1)\,dy\Pi(du)=
$$
$$=\frac{1}{A}\int_{|uz|>\beta}\left(\frac{\sin(e^{\,tA}|uz|)-\sin(|uz|)}{|uz|}-
(e^{\,tA}-1)\right)\Pi(du).
$$
Denote  $\gamma=\sup_{|y|>{(e^{tA}-1)\beta\over 2}}\left|{\sin
y\over y}\right|<1$. Then, for $|x|>\beta$,
$$
\frac{\sin(e^{tA}x)-\sin
x}{x}=(e^{tA}-1){\sin\left({(e^{tA}-1)x\over 2}\right)\over
{e^{tA}-1\over 2}}\cos\left({(e^{tA}+1)x\over 2}\right)\leq \gamma
(e^{tA}-1).
$$
Therefore
$$\int_0^tI_2(s,z)\,ds\leq-C_2\Pi(\{|uz|>\beta\}),\hbox{ where
}C_2=\frac{1-\gamma}{A}(e^{tA}-1)>0.
$$
Denote $C_3=\min(C_1\beta^2,C_2)$, then the estimates for
$\int_0^tI_{1,2}(s,z)\,ds$ given above yield
$$
|\phi_{X(t)}(z)|\leq \left({\beta\over
|z|}\right)^{-C_3\rho\left({\beta\over |z|}\right)},\quad z\in\Re.
$$
The latter estimate, under condition (iii), yields  (\ref{22}) and
therefore (i). The theorem is proved.

\begin{zau}\label{z1} The implication (ii) $\Rightarrow$ (iii) can be amplified with the following
statement: if $\lim\inf_{\eps\to
0+}\left[\varepsilon^2\ln\frac{1}{\varepsilon}\right]^{-1}\int_{\Re}
(u^2\wedge \eps^2)\Pi(du)=0,$ then for every  $t>0,p>1$ the
variable $X(t)$ does not possess a distribution density from the
class $L_p(\Re)$. In order to prove this fact, one should take
$\alpha={1\over 2}+{1\over 2p}\in(0,1)$ and a sequence  $\eps_n$
such that $\rho(\eps_n)\to 0$. Then the estimates analogous to
those given before provide that
$$
{P(X(t)\in(x_n,y_n))\over (y_n-x_n)^{2-2\alpha}}\to +\infty
$$
for $x_n=M(t,\eps_n)-\eps^\alpha_n,y_n=M(t,\eps_n)+\eps^\alpha_n$.
The latter convergence, together with H\"older inequality,
demonstrates that  $X(t)$ can not possess a distribution density
from the $L_{1\over 2\alpha -1}(\Re)=L_p(\Re).$
\end{zau}

Condition (iii) looks similar to the Kallenberg condition, but the
following example shows that these two conditions are remarkably
different.

\begin{ex} Let $\Pi=\sum_{n\geq 1}n\delta_{{1\over
n!}}$. Then
$$
\lim\inf_{\eps\to 0+}\rho(\eps)\geq \lim\inf_{\eps\to
0+}\left\{\Bigl[\ln{1\over \varepsilon}\Bigr]^{-1}
\Pi(|u|>\eps)\right\}\geq \lim\inf_{N\to+\infty}{1\over \ln N!}
\sum_{n\leq N-1} n \geq
$$
$$ \geq \lim\inf_{N\to+\infty}{N(N-1)\over 2 N\ln N}=+\infty,
$$
and condition (iii) holds true. One can check that the Kallenberg
condition fails, moreover,  we will show that the law of $Z(t)$ is
singular for every $t$. In order to do this, it is sufficient to
prove that $Ee^{izZ(t)}\not \to 0, z\to\infty$. But
$$
\lim_{N\to+\infty} \Bigl|Ee^{i 2\pi N!
Z(t)}\Bigl|=\lim_{N\to+\infty}\prod_{n>N}\Bigl|\exp\{tn(e^{i2\pi
N!\over n!}-1-{i2\pi N!\over n!})\}\Bigr|=1,
$$
that proves the needed statement. Thus, we have the following
interesting effect: the laws of  $Z(t)$ are singular, but the laws
of the values of the solution to (\ref{001}) with non-degenerated
drift  ($A\not=0$) possesses distribution densities of the class
$C_b^\infty$. One can say that the process $Z$ possess some
"hidden smoothness", that does not effect to the law of the
process itself, but becomes visible when this process is used as a
noise in an equation with a non-degenerate drift. Such an effect
is possible due to the difference between the Kallenberg condition
and (iii).
\end{ex}

From Proposition 2 and Theorem 1, one can make the general
conclusion that the conditions for \emph{existence} of the
distribution density for $X(t)$ on the one hand, and for
\emph{smoothness} of this density on the other, are essentially
different. This difference is well demonstrated by the following
example.

\begin{ex} Let $\Pi=\sum_{n\geq 1}\delta_{{1\over
n!}}$. Then $\Pi(\Re)=+\infty$, but
$$
\rho(\eps)\to 0, \quad \eps\to 0.
$$
For $A\not=0$, the solution to the equation (\ref{001}) possess
the distribution density, but this density is extremely irregular
in a sense that it does not belong to any $L_p(\Re),p>1$. These
two facts follow from Proposition 2 and Remark 1, respectively.
Remark that, in this example, the law of $Z$ is singular: one can
show this like it was done in Example 2. Thus, the current example
demonstrates one more version of a "regularization"\phantom{}
effect for the L\'evy process under the SDE with a non-degenerate
drift coefficient.
\end{ex}

\section{Multidimensional equation}

Let us introduce an auxiliary construction. Let a $\sigma$-finite
measure $\Pi$ to be defined on $\mathfrak{B}(\Re^d)$ with some
$d\in\mathbb{N}$. Consider the family  $\mathfrak{L}_\Pi=\{L$ is a
linear subspace of $\Re^d, \Pi(\Re^d\backslash L)<+\infty\}.$ It
is clear that if $L_1,L_2\in \mathfrak{L}_\Pi$ then $L_1\cap
L_2\in \mathfrak{L}_\Pi$. This yields that there exists a subspace
$L_\Pi\in \mathfrak{L}_\Pi$ such that $L_\Pi\subset L$ for every
$L\in\mathfrak{L}_\Pi$.

\begin{ozn} The subspace  $L_\Pi$ is called an essential linear support of the measure $\Pi.$
The measure $\Pi$ is said to be essentially linearly
non-degenerated if $L_\Pi=\Re^d$.
\end{ozn}

Remark that the condition on the measure  $\Pi$ to be essentially
linearly non-degenerated was imposed first in the paper
\cite{yamazato}, thus often it is called the Yamazato condition.

In this section, we study the local properties of the law of the
solution to the equation of the type \be\label{31} X(t)=\int_0^t
AX(s)\, ds+BW(t)+DZ(t),\quad t\geq 0 \ee  with  $A,B,D$ being an
$m\times m$-,$m\times k$- and $m\times d$-matrices respectively,
$W$ being a Wiener process in $\Re^k$ and the process $Z$ having
the form
$$
Z(t)=\int_0^t\int_{\|u\|_{\Re^d}>
1}u\nu(ds,du)+\int_0^t\int_{\|u\|_{\Re^d}\leq 1}u\tilde
\nu(ds,du),\quad t\geq 0
$$
with the L\'evy measure $\Pi$ being essentially linearly
non-degenerated. We have already seen that the solution $X$ to
equation (\ref{001}) linearly depends on $Z$: if $Z=Z_1+Z_2$ then
$X=X_1+X_2$, where $X_{1,2}$ denote the solutions to SDE
(\ref{001}) with $Z$ replaced by $Z_{1,2}$. When, moreover, $Z_1$
and $Z_2$ are independent and the process  $Z_2$ is the L\'evy
process without a diffusion component with its L\'evy measure
$\Pi_2$ being finite, then the law of $X_2(t)$ has an atom and
thus existence or smoothness of the distribution density for
$X(t)$ are equivalent to existence or smoothness of the
distribution density for $X_1(t)$. This allows one to remove from
the process $Z$ the part that is inessential in a sense of
Definition 1. Namely, one can  put $\nu_1$ equal to the
restriction of the point measure $\nu$ to  $\Re^+\times
(\Re^m\backslash L_\Pi)$ and define $Z_1$ by (\ref{0011}) with
$a=0, B=0$ and  $\nu$ replaced by $\nu_1$. One can easily see that
equation (\ref{001}) with $Z$ replaced by  $Z_{1}$ has the form
(\ref{31}) with $k=m, d=\mathrm{dim}\,L_\Pi$. Thus, equation
(\ref{001}) can be reduced to the form (\ref{31}), with $Z$
satisfying additionally the Yamazato condition.

If, in the equation (\ref{31}), $D=0$ then the following well
known Kalman controllability condition is necessary and sufficient
for the law of  $X(t),t>0$ to possess a smooth density (e.g.,
\cite{prat1}):
$$
 \quad \mbox{Rank}[B,AB,\ldots,A^{m-1}B]=m.
$$
Here $[B,AB,\ldots,A^{m-1}B]$ denotes an $m\times mk$-matrix of a
block type, composed from  the matrices з  $B,\ldots,A^{m-1}B$.
For the equation (\ref{31}), we write the analogous condition
$$
(\mathbf{H1}) \quad
\mbox{Rank}[B,AB,\ldots,A^{m-1}B,D,AD,\ldots,A^{m-1}D]=m,
$$
where $[B,AB,\ldots,A^{m-1}B,D,AD,\ldots,A^{m-1}D]$ is an $m\times
m(k+d)$-matrix composed form the matrices $B,\ldots,A^{m-1}B
,D,AD,\ldots,A^{m-1}D.$

Below, we denote  $S^d=\{l\in\Re^d, \|l\|_{\Re^d}=1\}$ (the unit
sphere in $\Re^d$). We introduce the multidimensional analogue of
the Kallenberg condition: \be\label{32}
\left[\varepsilon^2\ln\frac{1}{\varepsilon}\right]^{-1}\inf_{l\in
S^d}\int_{|(u,l)_{\Re^d}|\leq \eps} (u,l)^2_{\Re^d}\Pi(du)\to
+\infty, \quad \eps\to 0+.\ee We remark that this condition is a
new one.

\begin{thm} Let the L\'evy process $Z$ satisfy (\ref{32}). Then condition  $(\mathbf{H1})$
is sufficient for the law $X(t), t>0$ to possess a density from
the class $C_b^\infty$.
\end{thm}

\demo Like in the proof of Theorem 1, we will verify that the
characteristic function of $X(t)$ satisfies condition (\ref{22}).
We suppose that $\Pi(\|u\|_{\Re^d}>1)=0$, this obviously does not
restrict generality. The value of $X(t)$ is given as a sum of the
(independent) integrals over the Wiener process and the
compensated Poisson point measure. Thus the characteristic
function of $X(t)$ has the following explicit representation:
\begin{equation}\label{33}
\phi_{X(t)}(z)=\exp\left\{\int_0^t\left(-{1\over
2}\|B^*e^{(t-s)A^*}z\|_{\Re^k}^2+\int_{\Re^d}
\left[\exp\{i(e^{(t-s)A}Du,z)_{\Re^m}\}-1\right.\right.\right.-
$$
$$
\left.\left.\left.-i(e^{(t-s)A}Du,z)_{\Re^m}\right]
\Pi(du)\right)ds\right\},\quad z\in\Re^m,
\end{equation}
here $Q^*$ denotes adjoint matrix to $Q$ ($Q=A,B,\dots$). Then
\begin{equation} \label{003}
|\phi_{X(t)}(z)|=\exp{\left\{\int_0^t\left(-{1\over
2}\|B^*e^{(t-s)A^*}z\|_{\Re^k}^2+\int_{\Re^d}
\left[\cos{(e^{(t-s)A}Du,z)_{\Re^m}}-1\right]\Pi(du)\right)ds\right\}}.
\end{equation}
Denote $B(s,z)=B^*e^{\,sA^*}z,D(s,z)=D^*e^{\,sA^*}z$. We restrict
the domain of integration w.r.t. $u$ by the set
$\{|(D(s,z),u)_{\Re^d}|\leq 1\}$ and use inequality $1-\cos x \geq
Cx^2, |x|\leq 1, C=1-\cos 1>{1\over 2}$. We get \be\label{004}
|\phi_{X(t)}(z)|\leq\exp{\left\{-{1\over
2}\int_0^t\left(\|B(s,z)\|_{\Re^k}^2+\int_{|(D(s,z),u)_{\Re^d}|\leq
1}(D(s,z),u)_{\Re^d}^2\Pi(du)\right)ds\right\}}. \ee Denote
$$\Phi(r)=r^2\inf\limits_{l\in
S^m}\int_{|(u,l)_{\Re^d}|\leq{1\over
r}}(u,l)^2_{\Re^d}\Pi(du),\quad r>0,$$ remark that condition
(\ref{32}) is equivalent to the convergence  $\frac{\Phi(r)}{\ln
r}\to +\infty,\, r\to +\infty$. This notation allows us to rewrite
(\ref{004}) to the form
\begin{eqnarray} \label{005}
|\phi_{X(t)}(z)|\leq\exp{\left\{-{1\over
2}\int_0^t\left(\|B(s,z)\|_{\Re^k}^2+\Phi(\|D(s,z)\|_{\Re^d})\right)ds\right\}}.
\end{eqnarray}

\begin{lem} Under condition $(\mathbf{H1})$, for every given $t>0$ there exist
 $\alpha,\beta,\gamma>0$ such that
$$ \forall \, l\in S^m \quad \lambda\left\{0\leq s\leq t:\|B(s,l)\|_{\Re^k}>\alpha
\hbox{ або } \|D(s,l)\|_{\Re^d}>\beta\right\}\geq\gamma,$$ where
$\lambda$ is the Lebesgue measure on  $\Re$.
\end{lem}

\demo Suppose that the statement of the lemma does not hold true.
Then there exists a sequence $l_n\in S^m, n\geq 1$ such that
$$\lambda\left\{0\leq s\leq t:\|B(s,l_n)\|_{\Re^k}>\frac1n
\hbox{ or  } \|D(s,l_n)\|_{\Re^d}>\frac1n\right\}<\frac1n,\quad
n\geq1,$$ that means that both the sequences
$\{\|B(\cdot,l_n)\|_{\Re^k}\}$, $\{\|D(\cdot,l_n)\|_{\Re^d}\}$
converge in Lebesgue measure to the identical zero. Since $S^m$ is
a compact set, without a restriction of generality one can suppose
that $l_n\to l\in S^m$. But, for every $s\in [0,t]$, the functions
$B(s,\cdot), D(s,\cdot)$ are a linear and continuous ones, thus
the functions $\|B(\cdot,l)\|, \|D(\cdot,l)\|$ equal zero
$\lambda$-almost surely. Clearly, these functions are continuous,
and thus \be\label{35} B^*e^{sA^*}l=0,\quad D^*e^{sA^*}l=0,\quad
s\in[0,t]. \ee By taking the derivatives of (\ref{35}) by $s$ up
to the order $m-1$ and considering the values of the functions
$B(s,l),D(s,l)$ together with their derivatives  for  $s=0$, we
get
$$
B^*l=B^*A^*l=\dots=B^*(A^*)^{m-1}l=0, \quad
D^*l=D^*A^*l=\dots=D^*(A^*)^{m-1}l=0.
$$
The latter equality means that the rows of the matrix
$$[B,AB,\ldots,A^{m-1}B,D,AD,\ldots,A^{m-1}D]$$ are linearly dependent,
 with the coefficients of the dependence being equal to the coordinates of the vector $l$.
 This contradicts to condition $(\mathbf{H1})$. The lemma is
 proved.

 Now, we can complete the proof of Theorem 2. For a given
 $z\in\Re^m$, we put $l(z)={z\over \|z\|_\Re^m}$. Then
 $$
\lambda\left\{0\leq s\leq t:\|B(s,z)\|_{\Re^k}>\alpha
\|z\|_{\Re^m} \hbox{ або }
\|D(s,z)\|_{\Re^d}>\beta\|z\|_{\Re^m}\right\}= $$
$$
=\lambda\left\{0\leq s\leq t:\|B(s,l(z))\|_{\Re^k}>\alpha \hbox{
або } \|D(s,l(z))\|_{\Re^d}>\beta\right\}\geq\gamma.
 $$
The latter inequality and (\ref{005}) yield the estimate
\be\label{36} |\phi_{X(t)}(z)|\leq\exp{\left\{-{\gamma\over
2}\min\Big(\alpha\|z\|^2_{\Re^m},\Phi(\beta\|z\|^2_{\Re^m})\Big)\right\}},
\ee that, together with (\ref{32}), guarantees (\ref{22}). The
theorem is proved.

\begin{zau}
One can extend the result of Theorem 2 and describe in a more
details the asymptotic behavior of the derivatives of the density
$p_{X(t)}$ for $\|x\|_{\Re^m}\to\infty$. In order to make our
exposition transparent, we postpone the  discussion of this topic
to Section 5 below.
\end{zau}

\begin{zau}
In \cite{zab}, the statement is given (Theorem 3.1), being
analogous to Theorem 2. However, conditions imposed on the L\'evy
measure there  (Hypothesis 3.1) are somewhat superfluous and less
precise  than the multidimensional analogue (\ref{32}) of the
Kallenberg condition, used in the current paper.
\end{zau}

In \cite{zab}, Theorem 1.1,  it is proved that condition
$(\mathbf{H1})$ is sufficient for the law of the solution to
(\ref{31}) to be absolutely continuous, as soon  as  the jump
noise satisfies a multidimensional analogue of Sato condition (in
\cite{zab}, the case $B=0$ is considered, only). This statement
and Theorem 2 of the current paper show that $(\mathbf{H1})$ can
be naturally interpreted as the condition on the coefficients of
the equation that provides "preservation" of smoothness contained
in the additive noise $(W,Z)$. On the other hand, this condition
is satisfied for $A=0,B=0,D=I_{\Re^m}, d=m$. In this case
$X(t)=Z(t)-Z(0)$. Therefore, it is clear that condition
$(\mathbf{H1})$ does not provide a "regularization" effect,
analogous to the one of the one-dimensional equations with
non-degenerated drift obtained in the previous section.

Such kind of an effect, at least at the part of existence of the
density, is guaranteed by the following condition:
$$
(\mathbf{H2}) \quad \mathrm{Rank}\,[AD,\ldots,A^{m}D]=m.
$$
Although this condition contains the matrix $D$ as well as the
matrix $A$, we interpret it as an analogue of the condition on the
drift coefficient to be non-degenerate. We remark that this
condition is a new one, also.

\begin{thm} The following statements are equivalent:

(i) condition $(\mathbf{H2})$ holds true;

(ii) for an arbitrary solution to equation (8) with the process
$Z$ satisfying Yamazato condition, the law of the random vector
$X(t)$ is absolutely continuous for every $t>0$.
\end{thm}

\demo Let us prove the implication  (i) $\Rightarrow$ (ii) under
supposition that $B=0,\Pi(\|u\|_{\Re^d}>1)=0$. It was already
shown that such a supposition does not restrict generality since
the solution depends on the noise linearly. We use the sufficient
condition for absolute continuity of the law of a solution to SDE
with jump noise, given in Theorem 1.1 \cite{Me_jumps_UMZH}. This
condition is based on the construction proposed in
\cite{Me_TViMc}. We remark that, in  \cite{Me_jumps_UMZH},
 a general class of (non-linear) SDE's with jump noise is investigated under a specific
 moment condition (1.1). This condition is used in
\cite{Me_TViMc},\cite{Me_jumps_UMZH} in the proof of the
differentiability of the variable $X(t)$ w.r.t. certain group of
transformations of the Poisson point measure. For the equations
with an additive noise, such a differentiability holds true
without a specific moment condition (see
\cite{Nou_sim},\cite{Me}). Thus, we can apply the results obtained
in \cite{Me_jumps_UMZH} to the solution to (\ref{31}), not
requiring the moment condition  (1.1) \cite{Me_jumps_UMZH} to hold
true.

Statement A of Theorem 1.1 \cite{Me_jumps_UMZH} is formulated in
the terms of a certain subspace generated by a sequence of vector
fields, associated with the initial equation. In the partial case
of a linear equation (\ref{31}), these fields are defined as
$$ \Delta(u)=ADu,\quad
\mathfrak{L}(u)=\mathrm{Span}\Big\{\Lambda^k\Delta(u),
k\in\mathbb{Z}_+\Big\},\quad u\in\Re^d,\quad \Lambda v\eqdef-Av.
$$
By statement A of Theorem 1.1 \cite{Me_jumps_UMZH}, if for every
$l\in S^m$ \be\label{37} \Pi\Big(u: l \hbox{ is not orthogonal to
} \mathfrak{L}(u)\Big)=+\infty, \ee then the law of the solution
to (\ref{31}) is absolutely continuous.

Under condition $(\mathbf{H2})$, for every $l\in S^m$ there exists
proper subspace $L_l\subset \Re^d$ such that
$$
u\not\in L_l\Rightarrow \exists \,k\in\{1,\dots,m\}: \quad A^kD
u\not\perp l.
$$
Then $$ \Pi\Big(u: l \hbox{ is not orthogonal to }
\mathfrak{L}(u)\Big)\geq \Pi\Big(\Re^d\backslash L_l\Big).
$$
This, together with the Yamazato condition, provides  (\ref{37}).
The implication  (i) $\Rightarrow$ (ii) is proved.

Now, let us prove the inverse implication (ii) $\Rightarrow$ (i).
We put $B=0$. Let us prove that there exists a non-zero vector
 $l\in\Re^m$ such that \be\label{38}
(X(t),l)_{\Re^m}=(Z(t)-Z(0),D^*l)_{\Re^d},\quad t\geq 0. \ee If
$D=0,$ then (\ref{38}) trivially holds for every $l\in\Re^m$.
Thus, we suppose further that $D\not=0$. Under this supposition,
$\mathrm{Ker}\, D^*$ is a proper subspace of $\Re^m$. If
$(\mathbf{H2})$ does not hold, then there exists a non-zero vector
$l\in\Re^m$ such that \be\label{39} D^*A^*l=\dots=D^*(A^*)^{m}l=0,
\ee that means that the vectors $A^*l,\dots,(A^*)^{m}l$ belong to
the subspace $\mathrm{Ker}\, D^*$. Since the dimension of this
subspace does not exceed  $m-1$, there exist $k\leq m, c_1,\dots,
c_{k-1}\in\Re$: \be\label{310}
(A^*)^{k}l=\sum_{j=1}^{k-1}c_j(A^*)^{j}l.\ee By multiplying both
sides of (\ref{310}) on $(A^*)^{m+1-k}$ from the left,  and taking
into account that $(A^*)^{m+1-k+j}\mathrm{Ker}\, D^*,j\leq k-1$,
we get $(A^*)^{m+1}l\in \mathrm{Ker}\,D^*.$ Repeating these
considerations, we obtain that $(A^*)^{n}l\in \mathrm{Ker}\,
D^*,n\in\NN$ and thus  $e^{(t-s)A^*}l-l\in \mathrm{Ker}\, D^*,
0\leq s\leq t$. Then
$$
(X(t),l)_{\Re^m}=\int_0^t\int_{\|u\|_{\Re^d}\leq
1}(u,D^*e^{(t-s)A^*}l)_{\Re^d}\tilde \nu(ds,du)=
$$
$$
=\int_{0^t}\int_{\|u\|_{\Re^d}\leq
1}(u,D^*l)_{\Re^d}\tilde\nu(ds,du),
$$
that proves (\ref{38}).

If $D^*l=0,$ then (\ref{38}) immediately guarantees singularity of
the law $X(t)$ for every $Z$. Let us consider the case
$D^*l\not=0.$ Take the orthonormal basis $e_1,\dots,e_d$ in
$\Re^d$ in such a way that $e_1$ has the same direction with
$D^*l.$ Denote by $\gamma(r),r\geq 0$ the point in $\Re^d$ with
its coordinates (in this basis) equal $r,r^2,\dots, r^d$. A
standard argument using Vandermonde determinant provide that the
curve $\{\gamma(r),r\in\ax\}$ has at most $d$ intersection points
with every hyperplane in $\Re^d$. Thus the measure
$\Pi=\sum_{k\in\NN}\delta_{\gamma({1\over k!})}$ satisfies
Yamazato condition. On the other hand, the law of the variable
$(X(t),e_1)_{\Re^m}$ coincides with the law of the variable $Z(t)$
constructed in Example 2, and therefore is singular. Thus, the law
of  $X(t)$ is singular, also. The theorem is proved.

\begin{zau} Condition $(\mathbf{H1})$ involves non-trivially all the matrices
$A,B,D$, and thus the smoothness of the distribution density of
the solution to (\ref{31}) is provided by the diffusion and jump
noise conjointly. This differs from the statement (ii) of Theorem
3, where both the process $Z$ and the matrix $B$ are arbitrary.
The analogue of the condition $(\mathbf{H2})$ for the equation
with the fixed matrix  $B$ has the form
$$ (\mathbf{H2'}) \quad
\mathrm{Rank}\,[B,\dots, A^{m-1}B,AD,\ldots,A^{m}D]=m.
$$
The proof of necessity for this condition is totally analogous to
the proof  of necessity for $(\mathbf{H2})$, given in Theorem 3.
We can not use the results from the papers
\cite{Me_TViMc},\cite{Me_jumps_UMZH} in order to prove sufficiency
of this condition, since the equations with a diffusion component
are not considered in these papers. For the \emph{linear} SDE, the
approach developed in \cite{Me_TViMc},\cite{Me_jumps_UMZH},
without an essential changes, can be extended to SDE's with a
diffusion component, and using this approach one can prove
sufficiency of $(\mathbf{H2'})$. However, we do not give a
detailed exposition of this proof here, since, for such an
exposition, we would need to repeat the noticeable part of
\cite{Me_TViMc},\cite{Me_jumps_UMZH}.
\end{zau}
\begin{zau} It is well known that Kalman controllability condition
is, in fact, a version of the  H\"ormander hypoellipticity
condition, formulated for a separate class of linear diffusions.
Our condition $(\mathbf{H2})$ also has such an interpretation: in
the proof of sufficiency part of Theorem 3, we refer to Theorem
1.1 \cite{Me_jumps_UMZH}. The conditions given in the latter
theorem can be naturally considered as an analogue of the
H\"ormander hypoellipticity condition for SDE's with a jump noise.
\end{zau}

\begin{ex} (\cite{Me_jumps_UMZH}, Example 1.1) Consider the system
of equations
$$
\begin{cases}dX_1(t)=X_1(t)\, dt+dZ(t)\\
dX_2(t)=X_1(t)\, dt.
\end{cases}
$$
When  $Z$ is replaced by $W$ in this system, one get the well
known Kolmogorov's example of a two-dimensional diffusion
possessing smooth distribution density and being generated by a
one-dimensional Brownian motion. Initial system has the form
(\ref{31}) with $m=2,d=1, B=0$,
$$A=\left(\begin{array}{cc}1&0\\1&0\end{array}\right),\quad
D=\left(\begin{array}{c}1\\0\end{array}\right),\quad [D,AD]=
\left(\begin{array}{cc}1&1\\1&0\end{array}\right),\quad [AD,A^2D]=
\left(\begin{array}{cc}1&1\\1&1\end{array}\right).
$$
Condition $(\mathbf{H1})$ holds, but condition  $(\mathbf{H2})$
does not hold true. Thus, in the Kolmogorov's example, the
"preservation of smoothness" takes place, but the "regularization
effect" does not come into play.

Let us modify the Kolmogorov's example and consider the system of
equations
$$
\begin{cases}dX_1(t)=X_2(t)\, dt+dZ(t)\\
dX_2(t)=X_1(t)\, dt.
\end{cases}
$$
This system has the form (\ref{31}) with $m=2,d=1, B=0$,
$$A=\left(\begin{array}{cc}0&1\\1&0\end{array}\right),\quad
D=\left(\begin{array}{c}1\\0\end{array}\right),\quad [AD,A^2D]=
\left(\begin{array}{cc}0&1\\1&0\end{array}\right).
$$
For this system, condition $(\mathbf{H2})$ holds true. Thus, for
every process $Z$ with an infinite L\'evy measure (i.e., for a
process that has an infinite number of the jumps on every time
interval), the law of $X(t)=(X_1(t),X_2(t))$ in $\Re^2$ is
absolutely continuous.
\end{ex}

\section{Asymptotic properties of the derivatives of the distribution densities}

Together with the question on existence and smoothness of the
distribution density $p_{X(t)}(x),$ $x\in\Re^m$, it is natural to
study the limit behavior of the derivatives of this density for
$\|x\|_{\Re^m}\to +\infty$. In \cite{zab}, Remark 3.1, the problem
of integrability of the derivatives of the density $p_{X(t)}$ is
formulated in the connection with the investigation of the
smoothing properties of the semigroup generated by the process
$X$. In this section, we give a more strong version of Theorem 2,
that solve this problem completely.

Below, we denote by $\Sf(\Re^m)$ the Schwarz space of  infinitely
differentiable functions  $f:\Re^m\to \Re$ such that every
derivative of the function $f$, as $\|x\|_{\Re^m}\to\infty$, tends
to 0 faster than $\|x\|_{\Re^m}^{-n}$ for every  $n$.

\begin{thm} Consider equation (\ref{31}). If conditions (\ref{32}) and $(\mathbf{H1})$ hold,
then, for every $j_1,\dots,j_r\in\{1,\dots,m\}, r\in\NN, t>0$,
\be\label{41} {\prt^r\over \prt x_{j_1}\dots  \prt
x_{j_m}}p_{X(t)}\in L_1(\Re^m). \ee If, additionally, the L\'evy
measure of the process $Z$ satisfies the condition
\be\label{42}\int_{\|u\|_{\Re^d}>1}\|u\|_{\Re^d}^n\Pi(du)<+\infty,
\quad n\in\NN,\ee then  $p_{X(t)}\in\Sf(\Re^m), t>0$.
\end{thm}

\demo Let us consider first the case where the L\'evy measure
satisfies (\ref{42}). The Fourier transform is a bijective mapping
$\Sf(\Re^m)\to \Sf(\Re^m)$ (e.g., \cite{vlad}, \S 6.1). Thus, for
the proof of Theorem, it is sufficient to prove that every
derivative  of the characteristic function $\phi_{X(t)}(z),
z\in\Re^m$ tends to 0, as  $\|z\|_{\Re^m}\to\infty$, faster than
$\|z\|_{\Re^m}^{-n}$ for every $n$. This function has the
representation, analogous to (\ref{33}):
$\phi_{X(t)}=\exp[\psi_{X(t)}]$,
$$
\psi_{X(t)}(z)=\int_0^t\left(-{1\over
2}\|B^*e^{(t-s)A^*}z\|_{\Re^k}^2+\int_{\|u\|_{\Re^d>1}}
\left[\exp\{i(e^{(t-s)A}Du,z)_{\Re^m}\}-1\right] \Pi(du)\right.-
$$
$$
-\left.\int_{\|u\|_{\Re^d}\leq 1}
\left[\exp\{i(e^{(t-s)A}Du,z)_{\Re^m}\}-1-
i(e^{(t-s)A}Du,z)_{\Re^m}\right] \Pi(du)\right)\, ds,\quad
z\in\Re^m.
$$
Thus, every derivative of the function $\phi_{X(t)}$ has the form
$R\cdot\phi_{X(t)}$, where $R$ is some polynomial of the
derivatives of the function $\psi_{X(t)}$. We have already proved
in Theorem 2 that, under conditions (\ref{32}) and
$(\mathbf{H1})$,
$$
\phi_{X(t)}(z)=o(\|z\|_{\Re^m}^{-n}), \quad \|z\|_{\Re^m}\to
\infty,\quad n\in \NN.
$$
 Thus, it is enough to verify that every derivative of the function
$\psi_{X(t)}$ has at most polynomial growth as
$\|z\|_{\Re^m}\to\infty$. We have
$$
\frac{\partial}{\partial z_j}\psi_{X(t)}(z)=-\int_0^t
(B^*e^{(t-s)A^*}z,B^*e^{(t-s)A^*}e_j)_{\Re^k}ds+$$
$$+\int_0^t\int_{\|u\|_{\Re^d}\leq 1}i(e^{(t-s)A}Du,e_j)_{\Re^m}
\left[\exp\{i(e^{(t-s)A}Du,z)_{\Re^m}\}-1\right]\Pi(du)ds
$$
$$
+\int_0^t\int_{\|u\|_{\Re^d}>1}i(e^{(t-s)A}Du,e_j)_{\Re^m}
\exp\{i(e^{(t-s)A}Du,z)_{\Re^m}\}\Pi(du)ds ,\quad j=1,\dots,m,$$
where $e_j$ is the $j$-th basis vector in у $\Re^m$. Taking into
account the inequality $|e^{iz}-1|\leq |z|$, we get
\be\label{43}\left|\frac{\partial}{\partial
z_j}\psi_{X(t)}(z)\right|\leq C_1\left(\|z\|_{\Re^m}+
\|z\|_{\Re^m}\int_{\|u\|_{\Re^d}\leq
1}\|u\|_{\Re^d}^2\Pi(du)+\int_{\|u\|_{\Re^d}>
1}\|u\|_{\Re^d}\Pi(du)\right). \ee Here and below,  $C_r,
r=1,2,\dots$ are some constants that depend on coefficients
$A,B,D$ and time moment $t$. Next,
$$\frac{\partial^2}{\partial z_{j_1} \partial z_{j_2}}\psi_{X(t)}(z)=
-\int_0^t
(B^*e^{(t-s)A^*}e_{j_1},B^*e^{(t-s)A^*}e_{j_2})_{\Re^k}ds+$$
$$+\int_0^t\int_{\Re^d}i^2(e^{(t-s)A}Du,e_{j_1})_{\Re^m}(e^{(t-s)A}Du,e_{j_2})_{\Re^m}
\exp\{i(e^{(t-s)A}Du,z)_{\Re^m}\}\Pi(du)ds,
$$
$j_{1},j_{2}\in\{1,\dots,m\}$. Thus
\be\label{44}\left|\frac{\partial^2}{\partial z_{j_1}
\partial z_{j_2}}\psi(z)_{X(t)}(z)\right|\leq C_2\left(1+
\int_{\Re^d}\|u\|_{\Re^d}^2\Pi(du)\right),\quad
j_{1},j_{2}\in\{1,\dots,m\}. \ee At last, the partial derivatives
of the order  $r\geq3$ have the form
$$\frac{\partial^r}{\partial z_{j_1}\dots \partial z_{j_r}}\psi_{X(t)}(z)=
\int_0^t\int_{\Re^d}i^r\prod_{l=1}^r(e^{(t-s)A}Du,e_{j_l})_{\Re^m}
\exp\{i(e^{(t-s)A}Du,z)_{\Re^m}\}\Pi(du)ds,$$
$j_{1},\dots,j_{r}\in\{1,\dots,m\}$, that implies the estimate \be
\label{45} \left|\frac{\partial^r}{\partial z_{j_1}\dots
\partial z_{j_r}}\psi_{X(t)}(z)\right|\leq
C_r\left(\int_{\|u\|_{\Re^d}\leq
1}\|u\|_{\Re^d}^2\Pi(du)+\int_{\|u\|_{\Re^d}>
1}\|u\|_{\Re^d}^r\Pi(du)\right).\ee

It follows from (\ref{43}) -- (\ref{45}) that the first
derivatives of $\psi_{X(t)}$ have at most linear grows, while all
the higher derivatives are even bounded. Thus $\phi_{X(t)}\in
\Sf(\Re^m)$ and therefore $p_{X(t)}\in \Sf(\Re^m)$.

Now, let us consider the general case. Write
$W=W_1+W_2,Z=Z_1+Z_2$, $W_2=0, Z_2(t)=\int_0^t\int_{\|u\|_{\Re^d}>
1}u\nu(ds,du),$ and denote $X_{1,2}$ the solutions to SDE of the
type (\ref{31}) with $W,Z$ replaced by $W_{1,2}, Z_{1,2}$,
respectively. Then the solution to (8) has the form $X=X_1+X_2$,
and $X_1,X_2$ are independent. The density  $p_{X(t)}$ is equal
$$
p_{X(t)}(x)=\int_{\Re^m}p_{X_1(t)}(x-y)\mu_{X_2(t)}(dy), \quad
x\in\Re^m,
$$
where $\mu_{X_2(t)}$ denotes the law of $X_2(t)$. We have already
proved that $p_{X_1(t)}\in\Sf(\Re^m)$.  Thus,
$$
{\prt^r\over \prt x_{j_1}\dots  \prt
x_{j_m}}p_{X(t)}=\int_{\Re^m}{\prt^r\over \prt x_{j_1}\dots  \prt
x_{j_m}}p_{X_1(t)}(\cdot -y)\mu_{X_2(t)}(dy),
$$
$j_{1},\dots,j_{r}\in\{1,\dots,m\}, r\in\NN, t>0$, and
$$
\left\|{\prt^r\over \prt x_{j_1}\dots  \prt
x_{j_m}}p_{X(t)}\right\|_{L_1(\Re^m)}\leq
\int_{\Re^m}\left\|{\prt^r\over \prt x_{j_1}\dots \prt
x_{j_m}}p_{X_1(t)}(\cdot
-y)\right\|_{L_1(\Re^m)}\mu_{X_2(t)}(dy)=
$$
$$
=\left\|{\prt^r\over \prt x_{j_1}\dots  \prt
x_{j_m}}p_{X_1(t)}\right\|_{L_1(\Re^m)}.
$$
Theorem is proved.

\begin{zau} If (\ref{42}) does not hold, then $E\|Z(t)\|_{\Re^d}^n=+\infty$ for
some $n\in\NN$; the typical example here is provided by the stable
process with the index  $\alpha\in(0,2)$. Taking  $d=m,A=0, B=0,
D=I_{\Re^m},$ we get  $E\|X(t)\|_{\Re^m}^n=+\infty$. Therefore,
the condition (\ref{42}) is, in fact, necessary for the
distribution density of the solution $X(t)$ to belong to the
Schwarz space $\Sf(\Re^m)$.
\end{zau}

\section*{Conclusions} In the paper conditions are established,
allowing one to separate several questions that arise naturally
when the local properties of the laws of the solutions to SDE's
with jump noise are studied. The questions on \emph{existence} and
\emph{smoothness} of the distribution density appear to be
essentially different. Smoothness of the density is closely
related to the conditions on the behavior of the L\'evy measure of
the noise in the vicinity of the point  $0$ (Kallenberg condition
and its analogue (\ref{32}), condition (iii) of Theorem 1).
Conditions, necessary or sufficient for the density to exist, in
general, are much weaker. Moreover, the case of the equation that
contains a non-degenerate drift coefficient, appears to differ
essentially from the general one. For the equations with a
non-degenerate drift, on the contrary to the general ones, the
\emph{criteria} for existence and smoothness of the distribution
densities are available. In addition, non-degeneracy of the drift
coefficient makes possible the "regularization" of the
distribution of the L\'evy noise.

\end{document}